\documentclass[preprint,11pt]{elsarticle}

\usepackage[english]{babel}
\usepackage{amsmath, amsfonts, amssymb}

\newcommand{\ma}{\textit{Mathematica$^{\small{\circledR}}$}}

\usepackage{amsmath,amssymb,latexsym, color, amsfonts,  graphics}
\usepackage{theorem}

\usepackage{graphicx}

\newtheorem{theo}{Theorem}
\newtheorem{lemma}{Lemma}

\newtheorem{prop}{Proposition}
\newtheorem{nota}{Remark}
{\theorembodyfont{\rmfamily} }

\begin{document}

\begin{frontmatter}

\title{Asymptotic behavior of varying discrete Jacobi--Sobolev  orthogonal polynomials}

\tnotetext[label1]{The authors JFMM and JJMB are partially supported  by Research Group FQM-0229 (belonging to Campus of International Excellence CEIMAR). The author JFMM is funded by PPI Universidad de Almer\'{i}a. The author FM is partially supported by Direcci\'{o}n General de Investigaci\'{o}n, Ministerio de Econom\'{i}a y Competitividad Innovaci\'{o}n of Spain, Grant MTM2012-36732-C03-01. The author JJMB is partially supported by Direcci\'{o}n General de Investigaci\'{o}n, Ministerio de Ciencia e Innovaci\'{o}n of Spain and European Regional Development Found, grants MTM2011-28952-C02-01 and MTM2014-53963-P, and Junta de Andaluc\'{i}a (excellence grant P11-FQM-7276).}

\author[label2]{Juan F. Ma\~{n}as--Ma\~{n}as}
\ead{jmm939@ual.es}
\author[label3]{Francisco Marcell\'{a}n\tnoteref{label1}}
\ead{pacomarc@ing.uc3m.es}
\author[label2,label4]{Juan J. Moreno--Balc\'{a}zar\tnoteref{label1}}
\ead{balcazar@ual.es}

\address[label2]{Departamento de Matem\'{a}ticas, Universidad de Almer\'{\i}a, Spain}
\address[label3]{Instituto de Ciencias Matem\'aticas (ICMAT) and Departamento de
Matem\'{a}ticas, Universidad Carlos III de Madrid, Spain.}
\address[label4]{Instituto Carlos I de F\'{\i}sica
Te\'{o}rica y Computacional, Spain}

\begin{abstract}
In this contribution we deal with  a varying discrete Sobolev inner product involving the Jacobi weight. Our aim is to study the asymptotic properties of the corresponding orthogonal polynomials and the behavior of their zeros. We are interested in Mehler--Heine type formulae because they describe the essential  differences from the point of view of the asymptotic behavior between these Sobolev orthogonal polynomials and the Jacobi ones. Moreover, this  asymptotic behavior provides an approximation of the zeros of the Sobolev polynomials in terms of the zeros of other well--known special functions. We generalize some results appeared  in the literature very recently.

\end{abstract}


\end{frontmatter}

\section{Introduction} \label{s-int}

 One of the aims of  this paper is the study of the asymptotic behavior of sequences of polynomials   $\{Q_n^{(\alpha,\beta,M_n)}\}_{n\geq0}$ orthogonal  with respect to the inner product

\begin{equation}\label{projs}
(f,g)_{S,n}=\int_{-1}^1f(x)g(x)(1-x)^{\alpha}(1+x)^{\beta}dx+M_nf^{(j)}(1)g^{(j)}(1),
\end{equation}
where $\alpha>-1$, $\beta>-1,$ and $j\geq0.$

We assume that $\{M_n\}_{n\geq 0}$ is a sequence of nonnegative real numbers satisfying  \begin{equation}\label{sucesion}\lim_{n\to\infty}M_nn^{\gamma}=M>0,\end{equation}
where  $\gamma$ is a fixed real number.  Notice that this assumption is not very restrictive since the sequence $\{M_n\}_{n\geq 0}$ can behave asymptotically like any real power of the monomial $n.$

The main motivation to study this type of inner product arises from the papers \cite{alfandrrez} and \cite{alfmorpere}. In \cite{alfandrrez} the authors work with a  measure supported on $[-1,1].$ However, in \cite{alfmorpere} the authors deal with measures supported on an unbounded interval. In both cases the authors consider measures with nonzero absolutely continuous part, i.e., they work with the so--called continuous Sobolev orthogonal polynomials. The main topic in those papers is how to balance the Sobolev inner product to equilibrate the influence of the two measures in the asymptotic behavior of the corresponding orthogonal polynomials. This inspires us to consider the discrete Sobolev inner product
   $$
   (f, g)_S=\int fg d\mu_0+M\int f^{(j)}g^{(j)}d\mu_{1}=\int fg d\mu_0+Mf^{(j)}(c)f^{(j)}(c), $$
   which is a perturbation of a standard inner product. Now, making  $M$ dependent on $n$ we can study the influence of the perturbation on the asymptotic behavior of the orthogonal polynomials.  The literature on discrete Sobolev (or Sobolev--type) orthogonal polynomials is very wide, so we refer  the interested readers on this topic to  survey \cite{maryu} and the references therein.

   From here,  in \cite{mamamo} the authors found the asymptotic behavior of a family of orthogonal polynomials with respect to a varying Sobolev inner product similar to (\ref{projs}),  involving
the Laguerre weight $w(x)=x^{\alpha}e^{-x},\alpha>-1.$ We remark that the techniques used in \cite{mamamo} are not useful in this case, and now we need to use more powerful techniques based on those considered in \cite{dls}. More recently, in \cite{dls2015} the same authors have even improved these techniques in such a way that they have obtained relevant results for the orthogonal polynomials with respect to a non--varying discrete Sobolev inner product being $\mu_0$ a general measure.

 Previously, in \cite{morenoKrall}  J. J. Moreno--Balc\'{a}zar  obtained some results in this direction but only for the case $j=0$. Again, the method used in that paper does not allow to tackle our problem.

We want to emphasize that our objective is to establish that the size of the sequence $\{M_n\}_{n\geq 0}$ has an essential influence on the
asymptotic behavior of the orthogonal polynomials with respect to (\ref{projs}), but this influence is only local, that is, around the
point where we have introduced the perturbation. In our case, this point is located at $x=1$. Furthermore, we prove that this influence depends on the size of the sequence $\{M_n\}_{n\geq 0}$ and its relation with the parameter $\alpha$ in the Jacobi weight and the order of the derivative  in (\ref{projs}). It is important to remark that  for a
sequence $\{M_n\}_{n\geq 0},$ we have a sequence of orthogonal polynomials for each $n$, so we have a square tableau $\{Q_k^{(\alpha,\beta,M_n)} \}_{k\geq0}$. Here, we deal with the diagonal of this tableau, i.e. $\{Q_n^{(\alpha,\beta,M_n)}\}_{n\geq0}=\{ Q_{0}^{(\alpha,\beta,M_{0})}(x),Q_{1}^{(\alpha,\beta,M_{1})}(x),\dots, Q_{i}^{(\alpha,\beta,M_{i})}(x), \dots \}$. At this point, in order to simplify the notation, we will denote $Q_n^{(\alpha,\beta,M_n)}(x)=Q_{n}(x).$

 A second aim of this paper is to establish a simple asymptotic relation between the zeros of the Sobolev polynomials which are orthogonal with respect to  (\ref{projs}) and the zeros of combinations of Bessel functions of the first kind. This relation is deduced as an immediate consequence of Mehler--Heine formulae (Theorem \ref{amh}) and they have a numerical interest since we provide an estimate of the zeros of these polynomials.

Since Jacobi  classical orthogonal polynomials are involved in the varying inner product (\ref{projs}), we recall some  of their basic properties.  Jacobi polynomials are orthogonal with respect to the standard inner product
 $$
 (f,g)=\int_{-1}^1f(x)g(x)(1-x)^{\alpha}(1+x)^{\beta}dx, \quad \alpha, \beta>-1.
 $$

 In the sequel, we will work with the sequence   $\{P_n^{(\alpha,\beta)}\}_{n\geq 0},$  $\alpha>-1$ and $\beta>-1,$  normalized by (see \cite[f. (4.1.1)]{sz})

\begin{equation}\label{valorj1}
P_n^{(\alpha,\beta)}(1)=\binom{n+\alpha}{n}=\frac{\Gamma(n+\alpha+1)}{\Gamma(n+1)\Gamma({\alpha+1})}.
\end{equation}
The derivatives of Jacobi polynomials satisfy (see, \cite[f. (4.21.7)]{sz})
\begin{equation}\label{derjk}
(P_n^{(\alpha,\beta)}(x))^{(k)}=\frac{1}{2^k}\frac{\Gamma(n+\alpha+\beta+k+1)}{\Gamma(n+\alpha+\beta+1)}P_{n-k}^{(\alpha+k,\beta+k)}(x),\quad  k\geq 0.
\end{equation}
Using (\ref{valorj1}) and (\ref{derjk}), we deduce

\begin{equation}\label{valorderjk}
(P_n^{(\alpha,\beta)}(1))^{(k)}=\frac{1}{2^k}\frac{\Gamma(n+\alpha+\beta+k+1)}{\Gamma(n+\alpha+\beta+1)}
\frac{\Gamma(n+\alpha+1)}{\Gamma(n-k+1)\Gamma({\alpha+k+1})},
\end{equation}

\noindent where $(P_n^{(\alpha,\beta)}(1))^{(k)}$ denotes the  $k$th derivative of $P_n^{(\alpha,\beta)}$ evaluated at $x=1.$

We also note that the  squared norm of a Jacobi polynomial is (see, \cite[f. (4.3.3)]{sz}):

\begin{equation}\label{norma}
||P_n^{(\alpha,\beta)}||^2=\frac{2^{\alpha+\beta+1}}{2n+\alpha+\beta+1}\frac{\Gamma(n+\alpha+1)\Gamma(n+\beta+1)}{\Gamma(n+1)\Gamma(n+\alpha+\beta+1)}.
\end{equation}
Finally, we will use the  Mehler--Heine formula for  classical Jacobi polynomials

\begin{theo} (\cite[Th. 8.1.1]{sz}) \label{mhj} Let $\alpha,\beta>-1$. Then,
\begin{equation*}
\lim_{n\to \infty}n^{-\alpha}P_n^{(\alpha,\beta)}\left(\cos\left(\frac{x}{n}\right)\right)=\lim_{n\to\infty}\frac{1}{n^{\alpha}}P_n^{(\alpha,\beta)}
\left(1-\frac{x^2}{2n^2}\right)=(x/2)^{-\alpha}J_{\alpha}(x),
\end{equation*}
uniformly on compact subsets of $\mathbb{C}$. Here $J_{\alpha}(x)$ denotes the  Bessel function of the first kind, i.e.,
\begin{equation*}
J_{\alpha}(x)=\sum_{k=0}^{\infty} \frac{(-1)^k}{k!\Gamma(k+\alpha +1)} \left( \frac{x}{2} \right)^{2k+\alpha}.
\end{equation*}
\end{theo}

We will also use  the following limit related to Stirling formula (see, for example, \cite[f. (5.11.13)]{askey})

\begin{equation}\label{stirling}
\lim_{n\to\infty}\frac{n^{b-a}\Gamma(n+a)}{\Gamma(n+b)}=1.
\end{equation}

We introduce the following notation: If $a_n$ and $b_n$ are two sequences of real numbers, then  $a_n\approx b_n$ means that the sequence $\frac{a_n}{b_n}$ converges to 1.

 The paper is organized as follows. In Section \ref{s-vjsop} we provide some properties of the varying Jacobi--Sobolev orthogonal polynomials which are essential to establish  the Mehler--Heine asymptotics for these polynomials in Section \ref{s-asyzer}. Furthermore, as a consequence of this asymptotic formula we deduce the asymptotic behavior of the corresponding zeros. Thus, as we have commented previously, we can see the influence of the parameter $\gamma$, related to the size of the sequence $\{M_n\}_{n\ge 0},$  on the location of these zeros. Finally, in Section \ref{s-ne} we illustrate the results obtained in Section \ref{s-asyzer} with some numerical experiments.

\section{Varying Jacobi--Sobolev Orthogonal Polynomials} \label{s-vjsop}

It is well known that the classical Jacobi orthogonal polynomials,  $\{P_i^{(\alpha,\beta)}\}_{i=0}^{n},$
 constitute  a basis of the linear space $\mathbb{P}_n[x]$ of polynomials with real coefficients and degree at most $n.$ Therefore,  the Jacobi-Sobolev orthogonal polynomial of degree $n$, $Q_n(x)$, can be expressed as

\begin{equation*}
Q_n(x)=P_n^{(\alpha,\beta)}(x)+\sum_{i=0}^{n-1}a_{n,i}P_i^{(\alpha,\beta)}(x).
\end{equation*}
 Then, using well-known algebraic tools (see, for example, \cite[Sect. 2]{maRon}) we can deduce

\begin{equation}\label{expdesarrollada}
Q_n(x)=P_n^{(\alpha,\beta)}(x)-\frac{M_n\left(P_n^{(\alpha,\beta)}(1)\right)^{(j)}}{1+M_nK_{n-1}^{(j,j)}(1,1)}K_{n-1}^{(j,0)}(1,x),
\end{equation}
\noindent with $$K_n^{(j,k)}(x,y)=\sum_{i=0}^{n}\frac{\left(P_i^{(\alpha,\beta)}(x)\right)^{(j)}\left(P_i^{(\alpha,\beta)}(y)\right)^{(k)}}
{||P_i^{(\alpha,\beta)}(x)||^2}.$$

Next, we give a technical result useful for our purposes, interesting in itself though.

\begin{lemma}\label{lema1}
Let $\{Q_{n}\}_{n\geq 0}$ be the sequence of orthogonal polynomials with respect to (\ref{projs}) and $0\leq k \leq n$, then
\begin{description}
  \item[a)] \begin{equation}\label{comq/p}
\lim_{n\to\infty}\frac{(Q_n)^{(k)}(1)}{\left(P_n^{(\alpha,\beta)}(1)\right)^{(k)}}=\left\{
                                                                                     \begin{array}{ll}
                                                                                       \frac{k-j}{\alpha+j+k+1}, &\textrm{if} \quad \gamma<2(\alpha+2j+1), \\
                                                                                       \theta_{\alpha,\beta,j,k}, &\textrm{if} \quad \gamma=2(\alpha+2j+1), \\
                                                                                       1, &\textrm{if} \quad \gamma>2(\alpha+2j+1),
                                                                                     \end{array}
                                                                                   \right.
\end{equation}
where
\begin{equation} \label{theta}\theta_{\alpha,\beta,j,k}=\frac{M(k-j)+\Gamma^2(\alpha+j+1)2^{\alpha+\beta+2j+1}(\alpha+2j+1)(\alpha+j+k+1)}
{(\alpha+j+k+1)\left(M+\Gamma^2(\alpha+j+1)2^{\alpha+\beta+2j+1}(\alpha+2j+1)\right)}.\end{equation}
  \item[b)] $(Q_n,Q_n)_{S,n}\approx||P_n^{(\alpha,\beta)}||^2.$
\end{description}
\end{lemma}

\noindent \textbf{Proof.}  Kernel polynomials related to classical families of orthogonal polynomials and their derivatives have been  widely studied in the literature. Thus, we can claim  that the following limit exists,
  \begin{equation} \label{limkn} \lim_{n\to\infty} \frac{K_{n-1}^{(j,k)}(1,1)}{n^{2\alpha+2j+2k+2}}\in\mathbb{R}.\end{equation}  It is very easy to check it by using  Stolz's criterion, (\ref{valorderjk}), (\ref{norma}), (\ref{stirling}) and the fact that
$$n^{2\alpha+2j+2k+2}-(n-1)^{2\alpha+2j+2k+2}\approx (2\alpha+2j+2k+2)n^{2\alpha+2j+2k+1}.$$
 Thus,
\begin{eqnarray*}
\lim_{n\to\infty} \frac{K_{n-1}^{(j,k)}(1,1)}{n^{2\alpha+2j+2k+2}}&=&\lim_{n\to\infty} \frac{K_{n-1}^{(j,k)}(1,1)-K_{n-2}^{(j,k)}(1,1)}{n^{2\alpha+2j+2k+2}-(n-1)^{2\alpha+2j+2k+2}}\\
&=&\lim_{n\to\infty}\frac{\left(P_{n-1}^{(\alpha,\beta)}(1)\right)^{(k)}\left(P_{n-1}^{(\alpha,\beta)}(1)\right)^{(j)}}
{||P_{n-1}^{(\alpha,\beta)}||^2(2\alpha+2j+2k+2)n^{2\alpha+2j+2k+1}}\\
&=&\lim_{n\to\infty}\frac{C_{j,k} \Gamma(n+\alpha+\beta+j)\Gamma(n+\alpha+\beta+k)\Gamma(n+\alpha)\Gamma(n)}
{\Gamma(n-j)\Gamma(n+\alpha+\beta)\Gamma(n-k)\Gamma(n+\beta)n^{2\alpha+2j+2k}}\\
&=&C_{j,k}\in\mathbb{R},
\end{eqnarray*}
where
$$C_{j,k}=\frac{1}{\Gamma(\alpha+j+1)\Gamma(\alpha+k+1)2^{\alpha+\beta+j+k+1}(\alpha+j+k+1)}.$$

We will now prove part a) of the lemma, by (\ref{expdesarrollada})

\begin{eqnarray*}
\lim_{n\to\infty}\frac{Q_n^{(k)}(1)}{\left(P_n^{(\alpha,\beta)}(1)\right)^{(k)}}&=&
\lim_{n\to\infty}\left(1-\frac{M_nK_{n-1}^{(j,k)}(1,1)}{1+M_nK_{n-1}^{(j,j)}(1,1)}
\frac{\left(P_n^{(\alpha,\beta)}(1)\right)^{(j)}}{\left(P_n^{(\alpha,\beta)}(1)\right)^{(k)}}\right)\\
&=&\lim_{n\to\infty}\left(1-\frac{M_nn^{\gamma}\frac{1}{2^j}\frac{\Gamma(n+\alpha+\beta+j+1)}{\Gamma(n-j+1)\Gamma(\alpha+j+1)}
\frac{K_{n-1}^{(j,k)}(1,1)}{n^{2\alpha+2j+2k+2}}n^{2\alpha+2j+2k+2-\gamma}}
{\frac{1}{2^k}\frac{\Gamma(n+\alpha+\beta+k+1)}{\Gamma(n-k+1)\Gamma(\alpha+k+1)}\left(1+M_nn^{\gamma}
\frac{K_{n-1}^{(j,j)}}{n^{2\alpha+4j+2}}n^{2\alpha+4j+2-\gamma}\right)}\right).
\end{eqnarray*}
To simplify the computations we  introduce the following notation
\begin{eqnarray*}
a_n&=&M_n n^{\gamma},  \quad \textrm{by (\ref{sucesion}) we have}  \quad \lim_{n\to\infty}a_n=M,\\
 b_{n,j,k}&=&\frac{K_{n-1}^{(j,k)}(1,1)}{n^{2\alpha+2j+2k+2}}, \quad \mathrm{by \,\, (\ref{limkn}) \,\, we \,\, have } \quad  \lim_{n\to\infty}b_{n,j,k}=C_{j,k}.
\end{eqnarray*}
Then, the above limit becomes
\begin{eqnarray*}
&&\lim_{n\to\infty}\left(1-\frac{a_n b_{n,j,k}\frac{1}{2^j}\frac{\Gamma(n+\alpha+\beta+j+1)}{\Gamma(n-j+1)\Gamma(\alpha+j+1)}n^{2\alpha+2j+2k+2-\gamma}}
{\frac{1}{2^k}\frac{\Gamma(n+\alpha+\beta+k+1)}{\Gamma(n-k+1)\Gamma(\alpha+k+1)}\left(1+a_n b_{n,j,j} n^{2\alpha+4j+2-\gamma}\right)}\right)\\
&=&1-\frac{2^{k-j}\Gamma(\alpha+k+1)}{\Gamma(\alpha+j+1)}\ \lim_{n\to\infty}a_n\ \lim_{n\to\infty}b_{n,j,k}\ \lim_{n\to\infty}\frac{n^{2\alpha+2j+2k+2-\gamma}}{n^{2k-2j}(1+a_nb_{n,j,j}n^{2\alpha+4j+2-\gamma})} \\
&=&1-\frac{2^{k-j}\Gamma(\alpha+k+1)C_{j,k}M}{\Gamma(\alpha+j+1)}\ \lim_{n\to\infty}\frac{n^{2\alpha+4j+2-\gamma}}{(1+a_nb_{n,j,j}n^{2\alpha+4j+2-\gamma})}\\
&=&1-\frac{2^{k-j}\Gamma(\alpha+k+1)C_{j,k}M}{\Gamma(\alpha+j+1)}\ \lim_{n\to\infty}\frac{1}{\frac{1}{n^{2\alpha+4j+2-\gamma}}+MC_{j,j}}.
\end{eqnarray*}
Therefore, it is necessary to distinguish three cases according to the value of the parameter $\gamma.$ The value of this limit is:

\textit{Case}  $\gamma>2(\alpha+2j+1).$
\begin{eqnarray*}
1-\frac{2^{k-j}\Gamma(\alpha+k+1)C_{j,k}M}{\Gamma(\alpha+j+1)}\ \lim_{n\to\infty}\frac{1}{\frac{1}{n^{2\alpha+4j+2-\gamma}}+MC_{j,j}}=1.
\end{eqnarray*}

 \textit{Case} $\gamma<2(\alpha+2j+1).$
\begin{eqnarray*}
&&1-\frac{2^{k-j}\Gamma(\alpha+k+1)C_{j,k}M}{\Gamma(\alpha+j+1)}\ \lim_{n\to\infty}\frac{1}{\frac{1}{n^{2\alpha+4j+2-\gamma}}+MC_{j,j}}\\
&&=1-\frac{2^{k-j}\Gamma(\alpha+k+1)C_{j,k}M}{\Gamma(\alpha+j+1)}\ \frac{1}{MC_{j,j}}=1-\frac{\alpha+2j+1}{\alpha+j+k+1}=\frac{k-j}{\alpha+j+k+1}.
\end{eqnarray*}

 \textit{Case} $\gamma=2(\alpha+2j+1).$
\begin{eqnarray*}
&&1-\frac{2^{k-j}\Gamma(\alpha+k+1)C_{j,k}M}{\Gamma(\alpha+j+1)}\ \frac{1}{1+MC_{j,j}}\\
&&=1-\frac{M(\alpha+2j+1)}{(\alpha+j+k+1)\left(M+\Gamma^2(\alpha+j+1)2^{\alpha+\beta+2j+1}(\alpha+2j+1)\right)}\\
&&=\frac{M(k-j)+\Gamma^2(\alpha+j+1)2^{\alpha+\beta+2j+1}(\alpha+2j+1)(\alpha+j+k+1)}
{(\alpha+j+k+1)\left(M+\Gamma^2(\alpha+j+1)2^{\alpha+\beta+2j+1}(\alpha+2j+1)\right)}\\
&&=\theta_{\alpha,\beta,j,k}.
\end{eqnarray*}

Thus, we have proved a).  Now,  we are going to prove b).  Using standard arguments for Sobolev orthogonal polynomials we can deduce
\begin{equation*}
(Q_n,Q_n)_{S,n}=||P_n^{(\alpha,\beta)}||^2+\frac{M_n\left(\left(P_n^{(\alpha,\beta)}(1)\right)^{(j)}\right)^2}{1+M_nK_{n-1}^{(j,j)}(1,1)}.
\end{equation*}
Then,
\begin{eqnarray*}
\lim_{n\to\infty}\frac{(Q_n,Q_n)_{S,n}}{||P_n^{(\alpha,\beta)}||^2}&=&\lim_{n\to\infty}
\left(1+\frac{\left(\left(P_n^{(\alpha,\beta)}(1)\right)^{(j)}\right)^2}{||P_n^{(\alpha,\beta)}||^2}
\frac{M_n}{1+M_nK_{n-1}^{(j,j)}(1,1)}\right).
\end{eqnarray*}
 To establish b) it is enough to prove that

$$\lim_{n\to\infty}\left(\frac{M_n\left(\left(P_n^{(\alpha,\beta)}(1)\right)^{(j)}\right)^2}
{||P_n^{(\alpha,\beta)}||^2\left(1+M_nK_{n-1}^{(j,j)}(1,1)\right)}\right)=0. $$
Indeed, from (\ref{valorderjk}) and  (\ref{norma}) this limit can be expressed as
\begin{eqnarray*}
&\lim_{n\to\infty}&\left(\frac{M_n\left(\left(P_n^{(\alpha,\beta)}(1)\right)^{(j)}\right)^2}
{||P_n^{(\alpha,\beta)}||^2\left(1+M_nK_{n-1}^{(j,j)}(1,1)\right)}\right)=\\
&\lim_{n\to\infty}&\left(\frac{M_n\frac{1}{2^{2j}}\frac{\Gamma(n+\alpha+\beta+j+1)}{\Gamma(n-j+1)
\Gamma(\alpha+j+1)}\frac{\Gamma(n+\alpha+\beta+j+1)}{\Gamma(n+\alpha+\beta+1)}\frac{n^{-4j-2\alpha-\beta+\gamma}\Gamma(n+\alpha+1)}
{\Gamma(n-j+1)\Gamma(\alpha+j+1)}\frac{1}{n^{-4j-2\alpha-\beta+\gamma}}}
{\frac{2^{\alpha+\beta+1}}{2n+\alpha+\beta+1}\frac{n^{-\beta}\Gamma(n+\beta+1)}{\Gamma(n+1)}\frac{1}{n^{-\beta}}
\left(1+M_nK_{n-1}^{(j,j)}(1,1)\frac{n^{2\alpha+4j+2-\gamma}}{n^{2\alpha+4j+2-\gamma}}\right)}\right).
\end{eqnarray*}
 Again, to simplify the computations we introduce some notation
\begin{eqnarray*}
a_n&=&M_nn^{\gamma}, \,  \textrm{by (\ref{sucesion}) we have}   \lim_{n\to\infty}a_n=M,\\
b_n&=&\frac{\Gamma^2(n+\alpha+\beta+j+1)\Gamma(n+\alpha+1)n^{-4j-2\alpha-\beta}}{\Gamma(n-j+1)\Gamma^2(n+\alpha+\beta+1)},
 \, \textrm{then by (\ref{stirling}) }   \lim_{n\to\infty}b_n=1,\\
c_n&=&\frac{\Gamma(n+\beta+1)n^{-\beta}}{\Gamma(n+1)},\,  \textrm{then by (\ref{stirling})}   \lim_{n\to\infty}c_n=1,\\
d_n&=& M_nn^{\gamma}\frac{K_{n-1}^{(j,j)}(1,1)}{n^{2\alpha+4j+2}},\, \textrm{then using  (\ref{sucesion}) and (\ref{limkn}) we get}  \lim_{n\to\infty}d_n=MC_{j,j},\\
E_{\alpha,j}&=&\frac{1}{2^{2j}}\frac{1}{\Gamma^2(\alpha+j+1)}.
\end{eqnarray*}
In this way, for every $\gamma,$ the above limit is
\begin{eqnarray*}
\lim_{n\to\infty}\frac{E_{\alpha,j}a_nb_nn^{4j+2\alpha+\beta-\gamma}}{c_n\frac{2^{\alpha+\beta+1}n^{\beta}}{2n+\alpha+\beta+1}
\left(1+d_nn^{4j+2\alpha+2-\gamma}\right)}=\lim_{n\to\infty}\frac{E_{\alpha,j}a_nb_n(2n+\alpha+\beta+1)}{2^{\alpha+\beta+1}n^{2}
\left(\frac{c_n}{n^{4j+2\alpha+2-\gamma}}+c_nd_n\right)}=0,
\end{eqnarray*}
and we have just proved b).$\quad \Box$

\begin{nota}
Notice that taking into account b) in the above lemma, a) holds true when we consider orthonormal polynomials.
\end{nota}

 To tackle  Mehler--Heine asymptotics we need to expand the Sobolev polynomials $Q_n$ adequately. The following result gives us this expansion. In a more general framework it has been established in \cite[Th. 1]{dls2015}. The idea is that the coefficients $b_i(n)$ in (\ref{jsexp2}) can be obtained as a solution of a homogeneous linear system  of $j + 1$ equations and $j + 2$ unknowns. In our concrete case, we can compute explicitly the entries of the corresponding coefficient matrix.

\begin{prop}
There exists a family of real numbers $\{b_i(n)\}_{i=0}^{j+1}$, not identically zero, such that the following connection formula holds
\begin{equation}\label{jsexp2}
Q_n(x)=\sum_{i=0}^{j+1}b_i(n)(1-x)^{i}P_{n-i}^{(\alpha+2i,\beta)}(x), \quad   n\ge j+1.
\end{equation}
\end{prop}

\begin{lemma}
Let $\{b_i(n)\}_{i=0}^{j+1}$ be the coefficients in (\ref{jsexp2}).  Then $$\lim_{n\to\infty}b_i(n)=b_i\in\mathbb{R},\qquad i\in\{0, 1, \dots, j+1\}.$$
\end{lemma}
\textbf{Proof.} We take the $k$th derivative in (\ref{jsexp2}) and we evaluate the corresponding expression at $x = 1,$
\begin{eqnarray*}
Q_n^{(k)}(x)&=&\sum_{i=0}^{j+1}b_i(n)\sum_{s=0}^k\binom{k}{s}\left((1-x)^{i}\right)^{(s)}\left(P_{n-i}^{(\alpha+2i,\beta)}(x)\right)^{(k-s)}\\
&=&\sum_{i=0}^{j+1}b_i(n)\sum_{s=0}^{\min\{i,k\}}\binom{k}{s}(-1)^s\frac{i!}{(i-s)!}(1-x)^{i-s}\left(P_{n-i}^{(\alpha+2i,\beta)}(x)\right)^{(k-s)}.
\end{eqnarray*}
Then,
$$
Q_n^{(k)}(1)=\sum_{i=0}^{k}b_i(n)\binom{k}{i}(-1)^{i}i!\left(P_{n-i}^{(\alpha+2i,\beta)}(1)\right)^{(k-i)}.
$$
From Lemma \ref{lema1},  $\lim_{n\to\infty}\frac{Q_n^{(k)}(1)}{\left(P_n^{(\alpha,\beta)}(1)\right)^{(k)}}$ exists and its value depends on the value of parameter $ \gamma$ related to the size of the sequence $\{M_n\}_{n\ge 0}$, so

\begin{equation}\label{bii}
\frac{Q_n^{(k)}(1)}{\left(P_n^{(\alpha,\beta)}(1)\right)^{(k)}}=\sum_{i=0}^{k}b_i(n)\binom{k}{i}(-1)^{i}i!A_i(k,n)
\end{equation}
 with $A_i(k,n)=\frac{\left(P_{n-i}^{(\alpha+2i,\beta)}(1)\right)^{(k-i)}}{\left(P_n^{(\alpha,\beta)}(1)\right)^{(k)}}.$\\ It only remains to prove that there exists $\lim_{n\to\infty}A_i(k,n)\in\mathbb{R}$ and, therefore the coefficients  $\{b_i(n)\}_{i=0}^{j+1}$ are convergent. Indeed
\begin{eqnarray*}
\lim_{n\to\infty}A_i(k,n)&=&\lim_{n\to\infty}\frac{\frac{1}{2^{k-i}}\frac{\Gamma(n-i+\alpha+2i+\beta+k-i+1)}{\Gamma(n-i+\alpha+2i+\beta+1)}\frac{\Gamma(n-i+\alpha+2i+1)}{\Gamma(n-i-k+i+1)\Gamma(\alpha+2i+k-i+1)}}
{\frac{1}{2^k}\frac{\Gamma(n+\alpha+\beta+k+1)}{\Gamma(n+\alpha+\beta+1)}\frac{\Gamma(n+\alpha+1)}{\Gamma(n-k+1)\Gamma(\alpha+k+1)}}\\
&=&A_i(k,\alpha),
\end{eqnarray*}
 where we denote $\displaystyle A_i(k,\alpha)=\frac{2^{i}\Gamma(\alpha+k+1)}{\Gamma(\alpha+i+k+1)}. \quad \Box$

\begin{nota}
Let $b_i=\lim_{n\to \infty}b_i(n)$ with $i\in\{0,1,\dots, j+1\}.$    (\ref{bii})   is a recursive algorithm to compute $b_i.$
\begin{itemize}
  \item Step 1. For $k=0$ we obtain $b_0$ in a straightforward way.
  \item Step 2. For $k=1$ we deduce the value of $b_1$ from (\ref{bii}) using step 1. Similarly, for $k\ge 2$ we apply  (\ref{bii}) in a recursive way.
\end{itemize}
\end{nota}

\section{Asymptotics and zeros of varying Jacobi--Sobolev} \label{s-asyzer}

We focus our attention on the analysis of Mehler-Heine formulas for these discrete Jacobi--Sobolev orthogonal polynomials  because we want to know  how the discrete part in the inner product (\ref{projs}) influences  the asymptotic behavior of the corresponding
orthogonal polynomials. Furthermore, we will prove that this influence is related to  the size of the sequence $\{M_n\}_{n\geq 0}.$

\begin{theo}\label{amh}
For the sequence  $\{Q_n\}_{n\geq 0}$ the following Mehler--Heine formula holds
\begin{equation} \label{mhts}
\lim_{n\to\infty}\frac{Q_n(\cos(x/n))}{n^{\alpha}}=\lim_{n\to\infty}\frac{Q_n\left(1-\frac{x^2}{2n^2}\right)}{n^{\alpha}}=\left\{
                                                                         \begin{array}{ll}
                                                                           \phi_{\alpha}(x), & \hbox{if $\quad \gamma>2(\alpha+2j+1)$,} \\
                                                                           \psi_{\alpha,j}(x), & \hbox{if $\quad \gamma=2(\alpha+2j+1)$,} \\
                                                                           \varphi_{\alpha,j}(x), & \hbox{if $\quad \gamma<2(\alpha+2j+1)$,}
                                                                         \end{array}
                                                                       \right.
\end{equation}
uniformly on compact subsets of $\mathbb{C},$ where
$$
\phi_{\alpha}(x)=\left(\frac{x}{2}\right)^{-\alpha}J_{\alpha}(x),
$$
$$
\psi_{\alpha,j}(x)=\sum_{i=0}^{j+1}b_i2^{i}\left(\frac{x}{2}\right)^{-\alpha}J_{\alpha+2i}(x),
$$
with
$$
b_i=(-1)^i\frac{\frac{M(i-j)-\Gamma^2(\alpha+j+1)2^{\alpha+\beta+2j+1}(\alpha+2j+1)(\alpha+j+i+1)}
{(\alpha+j+i+1)(M+\Gamma^2(\alpha+j+1)2^{\alpha+\beta+2j+1}(\alpha+2j+1))}-
\Gamma(\alpha+i+1)\sum_{k=0}^{i-1}b_k\binom{i}{k}\frac{(-1)^{k}k! 2^{k}}{\Gamma(\alpha+i+k+1)}}
{i!\frac{2^i\Gamma(\alpha+i+1)}{\Gamma(\alpha+2i+1)}},
$$
for $0\leq i\leq j+1$, and
$$
\varphi_{\alpha,j}(x)=\sum_{i=0}^{j+1}b_i2^{i}\left(\frac{x}{2}\right)^{-\alpha}J_{\alpha+2i}(x),
$$
where  the coefficients  $b_i$ are computed as
$$
b_i= (-1)^i\frac{\frac{i-j}{\alpha+j+i+1}-\Gamma(\alpha+i+1)\sum_{k=0}^{i-1}b_k\binom{i}{k}\frac{(-1)^{k}k!2^{k}}{\Gamma(\alpha+i+k+1)}}
{ i!\frac{2^i\Gamma(\alpha+i+1)}{\Gamma(\alpha+2i+1)}}, \quad 0\leq i\leq j+1.
$$
Notice that in last two cases the coefficient $b_0$ is computed using the corresponding formula assuming $\sum_{i=0}^{-1} =0.$
\end{theo}
\textbf{Proof.}  Scaling and taking limits in (\ref{jsexp2})
\begin{eqnarray*}
\lim_{n\to\infty}\frac{Q_n\left(1-\frac{x^2}{2n^2}\right)}{n^{\alpha}}&=&
\lim_{n\to\infty}\frac{\sum_{i=0}^{j+1}b_i(n)\left(1-\left(1-\frac{x^2}{2n^2}\right)\right)^{i}P_{n-i}^{(\alpha+2i,\beta)}\left(1-\frac{x^2}{2n^2}\right)}{n^{\alpha}}\\
&=&\sum_{i=0}^{j+1}\lim_{n\to\infty}b_i(n)\ \lim_{n\to\infty}
\frac{\left(1-\left(1-\frac{x^2}{2n^2}\right)\right)^{i}P_{n-i}^{(\alpha+2i,\beta)}\left(1-\frac{x^2}{2n^2}\right)}{n^{\alpha}}\\
&=&\sum_{i=0}^{j+1}b_i2^{i}\left(\frac{x}{2}\right)^{-\alpha}J_{\alpha+2i}(x),
\end{eqnarray*}
 uniformly on compact subsets of $\mathbb{C}.$ Notice that in the last inequality  we have used Theorem \ref{mhj} written in the following way
$$
\lim_{n\to \infty} \frac{\left(\frac{x^2}{2n^2}\right)^i P_{n-i}^{(\alpha+2i,\beta)}\left(1-\frac{x^2}{2n^2}\right)}{n^{\alpha}} =2^i \left( \frac{x}{2} \right)^{-\alpha}J_{\alpha+2i}(x),$$
uniformly on compact subsets of $\mathbb{C},$ where $i$ is a fixed nonnegative integer number.

 Now, we distinguish three cases according to the value of the parameter $\gamma.$
\begin{itemize}
  \item If $\gamma>2(\alpha+2j+1),$ we are going to prove that $b_0=1$ and $b_i=0$ if $i\in\{1,2,\dots,j+1\}.$
\end{itemize}
 We can compute  $b_i$ from (\ref{bii}). If $k=0,$ then
\begin{eqnarray*}
\frac{Q_n(1)}{P_n^{(\alpha,\beta)}(1)}&=&b_0(n)A_0(0,n),
\end{eqnarray*}
Using   Lemma \ref{lema1} and taking limits, we obtain $b_0=1.$
If $k=1,$ then according to Lemma \ref{lema1} we have
$$
\frac{Q_n^{(1)}(1)}{\left(P_n^{(\alpha,\beta)}(1)\right)^{(1)}}=b_0(n)A_0(1,n)-b_1(n)A_1(1,n).
$$
Taking limits,
$$
1=1-b_1A_1(1,\alpha),\textrm{then  }
b_1=0.
$$
 Applying a recursive procedure we get $b_i=0$ for $i\in\{1,2,\dots,j+1\}$. To illustrate this procedure we consider the case $k=j+1.$ Thus, we have $b_i=0$ for $i\in\{1,2,\dots,j\}. $ Then,
\begin{eqnarray*}
\frac{Q_n^{(j+1)}(1)}{\left(P_n^{(\alpha,\beta)}(1)\right)^{(j+1)}}&=&b_0(n)A_0(1,n)+\sum_{i=1}^jb_i(n)\binom{j+1}{i}(-1)^{i}i!A_i(j+1,n)\\
&+&b_{j+1}(n)(-1)^{j+1}(j+1)!A_{j+1}(j+1,n).
\end{eqnarray*}
Taking limits,
$$
1=1+b_{j+1}(-1)^{j+1}(j+1)!A_{j+1}(j+1,\alpha), \textrm{then } b_{j+1}=0.$$
\begin{itemize}
  \item Case $\gamma=2(\alpha+2j+1).$ From (\ref{bii}) and  $k=0,$ we have
\end{itemize}
$$\frac{Q_n(1)}{P_n^{(\alpha,\beta)}(1)}=b_0(n)A_0(0,n).$$
 Taking limits when $n$ tends to infinity in the above expression, we get

$$b_0=\frac{-j M-\Gamma^2(\alpha+j+1)2^{\alpha+\beta+2j+1}(\alpha+2j+1)(\alpha+j+1)}
{(\alpha+j+1)(M+\Gamma^2(\alpha+j+1)2^{\alpha+\beta+2j+1}(\alpha+2j+1))}.$$

For $i\ge 1, $  we use  Lemma \ref{lema1} again and take limits. Thus, we deduce the coefficients  $b_i$ in a recursive way from  (\ref{bii}).

\begin{itemize}
  \item Case $\gamma<2(\alpha+2j+1).$ We can tackle this case in the same way as the case $\gamma=2(\alpha+2j+1). \quad  \Box$
\end{itemize}

 Next, we are going to study the zeros of the polynomials $\{Q_n\}_{n\geq 0}$ orthogonal with respect to (\ref{projs}). The following result was established for the non-varying case within a more general framework by H. G. Meijer in \cite[Th. 4.1]{mei} (see also \cite[Lemma 2]{allore}). Actually, that proof can be written in the same way for the varying case, so we omit it.
\begin{prop}\label{zeros-q}
The polynomial $Q_{n}(x),$ $n\geq 1,$   has $n$ real and simple zeros and at most one of them is located outside the interval $[-1,1]$.
\end{prop}
We can give more information about the location of the zeros. The case $j=0$ was considered in \cite{morenoKrall}. We notice that in that case all the zeros are in the interval $(-1,1).$ Thus, next we will assume $j> 0$ and we will denote by $y_{n,1}>y_{n,2}>\dots>y_{n,n-1}>y_{n,n}$ the zeros of $Q_n(x).$
\begin{prop} \label{zeros-location}  For $n$ large enough and $j>0$, we have
\begin{itemize}
  \item If $\gamma>2(\alpha+2j+1),$ then all zeros of $Q_n(x)$ are located in $(-1,1).$
  \item If $\gamma<2(\alpha+2j+1)$, then $y_{n,1}>1.$
  \item If $\gamma=2(\alpha+2j+1)$, then $y_{n,1}>1$ if and only if $$M>\frac{2^{\alpha+\beta+2j+1}(\alpha+j+1)(\alpha+2j+1)\Gamma^2(\alpha+j+1)}{j}$$
\end{itemize}
\end{prop}
\textbf{Proof.} We distinguish three cases, but essentially we use Lemma \ref{lema1} a) with $k=0$, and the fact that the leading coefficient of $Q_n$ is positive.  Then,
\begin{itemize}
  \item If $\gamma>2(\alpha+2j+1)$, then by Lemma \ref{lema1}  $Q_n(1)>0 $ for $n$ large enough. Therefore, taking into account Proposition \ref{zeros-q}, all the zeros are located in  $(-1,1)$.
  \item If $\gamma<2(\alpha+2j+1),$ then $Q_n(1)<0$ for $n$ large enough,  which implies that there is a zero of $Q_n$ greater than 1 and by Proposition \ref{zeros-q} it is the only one.
  \item If $\gamma=2(\alpha+2j+1)$, then $y_{n,1}>1$ if and only if $Q_n(1)<0 $ for $n$ large enough, and this only happens if and only if $$M>\frac{2^{\alpha+\beta+2j+1}(\alpha+j+1)(\alpha+2j+1)\Gamma^2(\alpha+j+1)}{j}. \quad \Box $$
\end{itemize}
Now  we  deduce the asymptotic behavior of the zeros of  $Q_n(x)$.
\begin{prop} \label{asybeh-zeros}
Let $y_{n,1}>y_{n,2}>\dots>y_{n,n-1}>y_{n,n}$ be the zeros of $Q_n(x)$ and $\phi_{\alpha}(x), \ \varphi_{\alpha,j}(x),$ and $\psi_{\alpha,j}(x)$ the functions defined in Theorem \ref{amh}. We assume $j>0.$
\begin{enumerate}
  \item If $\gamma>2(\alpha+2j+1),$ then   $$\lim_{n\to\infty} n\sqrt{2(1-y_{n,i})}=j_{\alpha,i},\quad i\ge 1,$$
where $j_{\alpha,i}$ denotes the $i$th positive zero of the Bessel function of the first kind.
  \item If $\gamma<2(\alpha+2j+1),$ then
   $$ \lim_{n\to\infty} y_{n,1}=1, \quad \lim_{n\to\infty}n\sqrt{2(1-y_{n,i})}=s_{\alpha,i-1},\quad i\ge 2,$$ where $s_{\alpha,i}$ denotes the $i$th positive zero of the function $\varphi_{\alpha,j}(x).$
  \item If $\gamma=2(\alpha+2j+1),$  we have two cases: \begin{enumerate}
                                            \item If $\displaystyle{M\leq\frac{2^{\alpha+\beta+2j+1}(\alpha+j+1)(\alpha+2j+1)\Gamma^2(\alpha+j+1)}{j}}$, then $y_{n,1}\leq1,$ for $n$ large enough, and  $$\lim_{n\to\infty}n\sqrt{2(1-y_{n,i})}=t_{\alpha,i},\quad i\ge 1,$$ where $t_{\alpha,i}$ denotes the $i$th positive zero of the function $\psi_{\alpha,j}(x).$
                                            \item If $\displaystyle{M>\frac{2^{\alpha+\beta+2j+1}(\alpha+j+1)(\alpha+2j+1)\Gamma^2(\alpha+j+1)}{j}}$, then $$ \lim_{n\to\infty} y_{n,1}=1, \quad   \lim_{n\to\infty}n\sqrt{2(1-y_{n,i})}=t_{\alpha,i-1},\quad i\ge 2,$$ where $t_{\alpha,i}$ denotes the $i$th positive zero of the function $\psi_{\alpha,j}(x).$
                                          \end{enumerate}
                                          \end{enumerate}
\end{prop}
\textbf{Proof.} It follows from Theorem \ref{amh}, Proposition \ref{zeros-location}, and Hurwitz's Theorem (see \cite[Th. 1.91.3]{sz}).$\quad \Box$

To illustrate Theorem \ref{amh} we are going to recover the case $j=0$  obtained in \cite{morenoKrall}. In that paper the author uses monic polynomials, and here we are considering a different normalization, i.e. the leading coefficient  of $Q_n$ is $$\frac{\Gamma(2n+\alpha+\beta+1)}{2^n \Gamma(n+1)\Gamma(n+\alpha+\beta+1)}.$$
Therefore, it is necessary to do some easy computations. We  use the relations (see, \cite[f.10.6.1]{askey}, \cite[6.1.18]{abra})
\begin{equation}\label{conexionbessel2}
J_{\alpha}(x)-\frac{2(\alpha+1)}{x}J_{\alpha+1}(x)=-J_{\alpha+2}(x),
\end{equation}
 as well as
\begin{equation}\label{duplicagamma}
\Gamma(2x)=\frac{\Gamma(x)\Gamma(x+\frac{1}{2})}{2^{1-2x}\sqrt{\pi}}.
\end{equation}

First, using (\ref{stirling}) and (\ref{duplicagamma}) we get
\begin{eqnarray*}
\frac{\Gamma(2n+\alpha+\beta+1)}{2^n \Gamma(n+1)\Gamma(n+\alpha+\beta+1)}&\approx&
\frac{2^{n+\alpha+\beta}}{\sqrt{\pi}}\frac{\Gamma\left(n+\frac{\alpha}{2}+\frac{\beta}{2}+\frac{1}{2}\right)}{ \Gamma(n+1)}
\frac{\Gamma\left(n+\frac{\alpha}{2}+\frac{\beta}{2}+1\right)}{\Gamma(n+\alpha+\beta+1)}\\
&\approx&\frac{2^{n+\alpha+\beta}}{n^{\frac{1}{2}}\sqrt{\pi}}.
\end{eqnarray*}
In \cite{morenoKrall} it was obtained
\begin{equation*}
\lim_{n\to\infty}\frac{2^n \hat{P}_n^{(\alpha,\beta,M_n)}(\cos(x/n))}{n^{\alpha+1/2}}=\left\{
                                                                                        \begin{array}{ll}
                                                                                          -2^{-\beta}\sqrt{\pi}x^2z_{\alpha+2}(x), & \textrm{if }  \hbox{$\gamma<2\alpha+2$,} \\
                                                                                          -2^{-\beta}\sqrt{\pi}(z_{\alpha}(x)+a_{\alpha,\beta,M}z_{\alpha+1}(x)), & \textrm{if } \hbox{$\gamma=2\alpha+2$,} \\
                                                                                          2^{-\beta}\sqrt{\pi}z_{\alpha}(x), & \textrm{if } \hbox{$\gamma>2\alpha+2$,}
                                                                                        \end{array}
                                                                                      \right.
\end{equation*}
where
\begin{eqnarray*}
z_{\alpha}(x)&=&x^{-\alpha}J_{\alpha}(x),\\
a_{\alpha,\beta,M}&=&\frac{-2M(\alpha+1)}{M+2^{\alpha+\beta+1}\Gamma(\alpha+2)\Gamma(\alpha+1)},
\end{eqnarray*}
 and $\{\hat{P}_n^{(\alpha,\beta,M_n)}\}_{n\geq0}$ denotes the sequence of monic  polynomials which are orthogonal with respect to (\ref{projs}) with $j=0.$
This result can be written as follows
\begin{equation} \label{mhjuanjo}
\lim_{n\to\infty}\frac{2^{n+\alpha+\beta} \hat{P}_n^{(\alpha,\beta,M_n)}(\cos(x/n))}{n^{\alpha+1/2}\sqrt{\pi}}=\left\{
                                                                                        \begin{array}{ll}
                                                                                          -2^{\alpha}x^2z_{\alpha+2}(x), & \textrm{if }  \hbox{$\gamma<2\alpha+2$,} \\
                                                                                          -2^{\alpha}(z_{\alpha}(x)+a_{\alpha,\beta,M}z_{\alpha+1}(x)), & \textrm{if } \hbox{$\gamma=2\alpha+2$,} \\
                                                                                          2^{\alpha}z_{\alpha}(x), & \textrm{if } \hbox{$\gamma>2\alpha+2$.}
                                                                                        \end{array}
                                                                                      \right.
\end{equation}
 We can observe that
$$\frac{2^{n+\alpha+\beta} \hat{P}_n^{(\alpha,\beta,M_n)}(\cos(x/n))}{n^{\alpha+1/2}\sqrt{\pi}}\approx\frac{Q_n(\cos(x/n))}{n^{\alpha}}.$$

Therefore, it only remains to compare the limit functions in (\ref{mhjuanjo}) and (\ref{mhts}). The case $\gamma>2\alpha+2$ is trivial. We pay attention to the other two cases.

\begin{itemize}

  \item  $\gamma<2\alpha+2.$

In this case  $b_0=0$ and $b_1=-1/2.$ Thus we have

$$
\varphi_{\alpha,0}(x)=-\left(\frac{x}{2}\right)^{-\alpha}J_{\alpha+2}(x)=-2^{\alpha}x^2x^{-\alpha-2}J_{\alpha+2}(x)=
-2^{\alpha}x^2z_{\alpha+2}.$$
  \item  $\gamma=2\alpha+2.$

In this case,
\begin{eqnarray*}
b_0&=&-\frac{\Gamma^2 (\alpha+1)2^{\alpha+\beta+1} (\alpha+1)}{M+\Gamma^2(\alpha+1)2^{\alpha+\beta+1}(\alpha+1)},\\
b_1&=&\frac{M}{2(M+\Gamma^2(\alpha+1) 2^{\alpha+\beta+1}(\alpha+1))}.
\end{eqnarray*}
By using (\ref{conexionbessel2}) we deduce
\begin{eqnarray*}
 \psi_{\alpha,0}(x)&=&b_0\left(\frac{x}{2}\right)^{-\alpha}J_{\alpha}(x)+2b_1\left(\frac{x}{2}\right)^{-\alpha}J_{\alpha+2}(x)\\
&=&\frac{-\Gamma^2(\alpha+1)2^{\alpha+\beta+1}}{M+\Gamma^2(\alpha+1)2^{\alpha+\beta+1}(\alpha+1)}\left(\frac{x}{2}\right)^{-\alpha}J_{\alpha}(x)\\
&+&\frac{M}{M+\Gamma^2(\alpha+1)2^{\alpha+\beta+1}(\alpha+1)}\left(\frac{x}{2}\right)^{-\alpha}J_{\alpha+2}(x)\\
&=&-\left(\frac{x}{2}\right)^{-\alpha}J_{\alpha}(x)+\frac{M(\alpha+1)}{M+2^{\alpha+\beta+1}\Gamma^2(\alpha+1)(\alpha+1)}\left(\frac{x}{2}\right)^{-\alpha-1}J_{\alpha+1}(x)\\
&=&-2^{\alpha}(z_{\alpha}(x)+a_{\alpha,\beta,M}z_{\alpha+1}(x)).
\end{eqnarray*}
\end{itemize}

\section{Numerical Experiments} \label{s-ne}

In this section we illustrate the  previous results  on the zeros of the polynomials $Q_n$  with some  numerical experiments where we have taken $j=3 $ for all of them. Thus, we are dealing with the varying Sobolev inner product
$$
(f,g)_{S,n}=\int_{-1}^1f(x)g(x)(1-x)^{\alpha}(1+x)^{\beta}dx+M_nf^{(3)}(1)g^{(3)}(1).
$$

We have used  the mathematical software \ma\textit{8.0} for the computations. In all the numerical experiments we have computed
 the four largest zeros of the polynomials $Q_n(x)$ and the corresponding scaled zeros for several values of $n.$  We only show one example for each possible case. In the tables about the scaled zeros we show their asymptotic behavior such as it is described in Proposition \ref{asybeh-zeros}.

\begin{itemize}
  \item Case $\gamma>2(\alpha+2j+1).$
\end{itemize}
We choose the following values:
$$\alpha=3,\quad  \beta=1,\quad \gamma=25, \quad \mathrm{and} \quad M_n=\frac{3 e^{n}}{(6e^n+4)n^{\gamma}}.$$

 It was proved in Theorem \ref{amh} that in this case the Mehler--Heine formula for the polynomials $Q_n$ is the same one as for the classical Jacobi polynomials. This behavior is due to the negligible influence of the sequence of masses $\{M_n\}_{n\geq 0}$ on the asymptotics.  Obviously, as it was stated in Proposition \ref{asybeh-zeros}, this determines the asymptotic behavior of the zeros which is illustrated in Table \ref{tabla1} and Table \ref{tabla2}.

\begin{table}[h]
\begin{center}
\scalebox{0.80}{
\begin{tabular}{|c|c|c|c|c|}
\hline
  & $y_{n,4}$ & $y_{n,3}$ & $y_{n,2}$  & $y_{n,1}$ \\ \hline
 $n=150$ & $0.994346$ & $0.99636$ & $0.997952$ & $0.999125$            \\ \hline
 $n=250$ & $0.997937$ & $0.998672$ & $0.999254$ & $0.999681$              \\ \hline
 $n=500$ & $0.999479$ & $0.999665$ & $0.999811$ & $0.999919$                            \\ \hline
\end{tabular}}
\end{center}
\begin{center}
\caption[Table 1]{\textbf{Case $\gamma=25>2(\alpha+2j+1),$}\\}\label{tabla1}
 $j=3,\ \alpha=3,\ \beta=1,\ \gamma=25,\ M_n=\frac{3 e^{n}}{(6e^n+4)n^{\gamma}}.$
\end{center}
\end{table}

\begin{table}[h]
\begin{center}
\scalebox{0.80}{
\begin{tabular}{|c|c|c|c|c|}
\hline
  & $n\sqrt{2(1-y_{n,1})}$ & $n\sqrt{2(1-y_{n,2})}$ & $n\sqrt{2(1-y_{n,3})}$  & $n\sqrt{2(1-y_{n,4})}$ \\ \hline
 $n=150$ & $6.27524$ & $9.59956$ & $12.7982$ & $15.9503$            \\ \hline
 $n=250$ & $6.31687$ & $9.66386$ & $12.885$  & $16.0602$              \\ \hline
 $n=500$ & $6.34839$ & $9.71233$ & $12.9501$ & $16.1421$                            \\ \hline
\bf{Limit} &$\mathbf{j_{3,1}}=6.38016$ &$\mathbf{j_{3,2}}=9.76102 $ & $\mathbf{j_{3,3}} =13.0152 $ & $\mathbf{j_{3,4}} = 16.2235 $   \\ \hline
\end{tabular}}

\end{center}
\begin{center}
\caption[Table 2]{\textbf{Case $\gamma=25>2(\alpha+2j+1),$}\\ }\label{tabla2}
$j=3,\ \alpha=3,\ \beta=1,\ \gamma=25,\ M_n=\frac{3 e^{n}}{(6e^n+4)n^{\gamma}}.$
\end{center}
\end{table}

\begin{itemize}
  \item Case $\gamma<2(\alpha+2j+1).$
\end{itemize}

According to Theorem \ref{amh} the limit function in the Mehler--Heine formula is given by   $\varphi_{\alpha,3}(x)=\sum_{i=0}^{4}b_i2^{i}\left(\frac{x}{2}\right)^{-\alpha}J_{\alpha+2i}(x),$ where the coefficients $b_i, \ 0\leq  i\leq 4,$ can be computed from  Theorem \ref{mhj}.We choose the following values:
$$\alpha=3,\quad  \beta=1,\quad \gamma=4, \quad  \mathrm{and}  \quad M_n=\frac{7 \ln(n+1)+5}{(3+\ln(n^2))n^{\gamma}}.$$

In Table \ref{tabla3} we can see that the largest zero is greater than 1 for $n$ large enough according to Proposition \ref{zeros-location}.  Table \ref{tabla4}  shows the asymptotic behavior of the scaled zeros given in Proposition \ref{asybeh-zeros}.

\begin{table}[ht]
\begin{center}
\scalebox{0.80}{
\begin{tabular}{|c|c|c|c|c|}
\hline
  & $y_{n,4}$ & $y_{n,3}$ & $y_{n,2}$  & $y_{n,1}$ \\ \hline
 $n=150$ & $0.994574$ & $0.996593$ & $0.998169$ & $0.999286$            \\ \hline
 $n=250$ & $0.998176$ & $0.998915$ & $0.999497$ & $1.0016$              \\ \hline
 $n=500$ & $0.999554$ & $0.999739$ & $0.999883$ & $1.0014$                            \\ \hline
\end{tabular}}
\end{center}
\begin{center}
\caption[Tabla 3]{\textbf{Case $\gamma=4<2(\alpha+2j+1),$}\\ } \label{tabla3}
$j=3,\ \alpha=3,\ \beta=-1/2,\ \gamma=4,\ M_n=\frac{7 \ln(n+1)+5}{(3+\ln(n^2))n^{\gamma}}.$
\end{center}
\end{table}
\begin{table}[ht]
\begin{center}
\scalebox{0.80}{
\begin{tabular}{|c|c|c|c|}
\hline
  &  $n\sqrt{2(1-y_{n,2})}$ & $n\sqrt{2(1-y_{n,3})}$  & $n\sqrt{2(1-y_{n,4})}$ \\ \hline
 $n=150$ & $9.07735$ & $12.382$ & $15.6257$            \\ \hline
 $n=250$& $7.92964$ & $11.6463$ & $15.1011$               \\ \hline
 $n=500$& $7.6415$ & $11.4238$ & $14.9355$                             \\ \hline

\bf{Limit}  &$\mathbf{s_{3,1}}=7.64622$ &$\mathbf{s_{3,2}}=11.4432 $ & $\mathbf{s_{3,3}} =14.9699 $    \\ \hline
\end{tabular}}
\end{center}
\begin{center}
\caption[Tabla 4]{\textbf{Case $\gamma=4<2(\alpha+2j+1),$}\\ }\label{tabla4}
 $j=3,\ \alpha=3,\ \beta=-1/2,\ \gamma=4,\ M_n=\frac{7 \ln(n+1)+5}{(3+\ln(n^2))n^{\gamma}}.$
\end{center}
\end{table}

\begin{itemize}
  \item Case $\gamma=2(\alpha+2j+1).$ \end{itemize}

According to Theorem \ref{amh} the limit function in the Mehler--Heine formula is given by $\psi_{\alpha,3}(x)=\sum_{i=0}^{4}b_i2^{i}\left(\frac{x}{2}\right)^{-\alpha}J_{\alpha+2i}(x),$ where the coefficients $b_i, \ 0\leq  i\leq 4,$ can be computed again from  Theorem \ref{mhj}.
We choose the following values:
$$\alpha=\beta=-9/10, \quad \gamma=61/5=12.2, $$
and we denote by $V$ the quantity which appears in Proposition \ref{asybeh-zeros}, i.e.
$$V=\frac{2^{\alpha+\beta+2j+1}(\alpha+j+1)(\alpha+2j+1)\Gamma^2(\alpha+j+1)}{j}.$$
Thus, with this data
$$V=2^{1/5}\frac{15128}{75}\Gamma^2\left(\frac{31}{10}\right)\simeq 1119.0037947. $$
Now we take
$$
M_n=\frac{Mn^2(n-1/2)(n+2)}{n^{\gamma+4}}=\frac{Mn^2(n-1/2)(n+2)}{n^{81/5}}.$$
According to Proposition \ref{asybeh-zeros} we have two possible choices of $M$ which determine two different asymptotic behaviors of the zeros. In Table \ref{tabla5} and Table \ref{tabla6}  we show the case $M\le V$ where $M=5.$ We can see that the largest zero of $Q_n$ is always lesser than 1. However, when $M>V$ then $y_{n,1} >1$ for $n$ large enough and this is illustrated in Table \ref{tabla7} for $M=10^6$. In Table \ref{tabla8} the asymptotic behavior of the scaled zeros is shown.
\begin{table}[h]
\begin{center}
\scalebox{0.80}{
\begin{tabular}{|c|c|c|c|c|}
\hline
  & $y_{n,4}$ & $y_{n,3}$ & $y_{n,2}$  & $y_{n,1}$ \\ \hline
 $n=150$ & $0.99778$ & $0.99854$ & $0.999585$ & $0.999991$            \\ \hline
 $n=250$ & $0.999142$ & $0.999585$ & $0.999871$ & $0.999997$              \\ \hline
 $n=500$ & $0.999786$ & $0.99985$ & $0.999968$ & $0.999999$                            \\ \hline
\end{tabular}}
\end{center}
\begin{center}
\caption[Tabla 5]{\textbf{Case $\gamma=61/5=2(\alpha+2j+1),$}\, $M=5\le V$\\ } \label{tabla5}
 $j=3,\ \alpha=-9/10,\ \beta=-9/10,\ \gamma=61/5,\ M_n=\frac{5n^2(n-1/2)(n+2)}{n^{\gamma+4}}.$
\end{center}
\end{table}
\begin{table}[h]
\begin{center}
\scalebox{0.80}{
\begin{tabular}{|c|c|c|c|c|}
\hline
  & $n\sqrt{2(1-y_{n,1})}$ & $n\sqrt{2(1-y_{n,2})}$ & $n\sqrt{2(1-y_{n,3})}$  & $n\sqrt{2(1-y_{n,4})}$ \\ \hline
 $n=150$ & $0.649565$ & $4.02672$ & $7.20558$ & $10.3659$            \\ \hline
 $n=250$ & $0.64887$ & $4.02249$ & $7.19831$ & $10.3561$              \\ \hline
 $n=500$ & $0.64853$ & $4.01929$ & $7.19273$ & $10.3484$                            \\ \hline
\bf{Limit} &$\mathbf{t_{0,1}}=0.648561$ &$\mathbf{t_{0,2}}=4.01985 $ & $\mathbf{t_{0,3}} =7.19169 $ & $\mathbf{t_{0,4}} = 10.3446 $   \\ \hline
\end{tabular}}
\end{center}
\begin{center}
\caption[Tabla 6]{\textbf{Case $\gamma=61/5=2(\alpha+2j+1),$} \, $M=5\le V$\\ }\label{tabla6}
 $j=3,\ \alpha=-9/10,\ \beta=-9/10,\ \gamma=61/5,\ M_n=\frac{5n^2(n-1/2)(n+2)}{n^{\gamma+4}}.$
\end{center}
\end{table}
\begin{table}[h]
\begin{center}
\scalebox{0.80}{
\begin{tabular}{|c|c|c|c|c|}
\hline
  & $y_{n,4}$ & $y_{n,3}$ & $y_{n,2}$  & $y_{n,1}$ \\ \hline
 $n=150$ & $0.996412$ & $0.999306$ & $0.999978$ & $1.00042$            \\ \hline
 $n=250$ & $0.99931 $ & $0.999739$ & $0.999991$ & $1.00009$              \\ \hline
 $n=500$ & $0.999818$ & $0.999928$ & $0.999999$ & $1.000001$                            \\ \hline
\end{tabular}}
\end{center}
\begin{center}
\caption[Tabla 7]{\textbf{Case $\gamma=61/5=2(\alpha+2j+1),$} $M=10^6>V$\\ } \label{tabla7}
 $j=3,\ \alpha=-9/10,\ \beta=-9/10,\ \gamma=61/5,\ M_n=\frac{10^6n^2(n-1/2)(n+2)}{n^{\gamma+4}}.$
\end{center}
\end{table}
\begin{table}[h]
\begin{center}
\scalebox{0.80}{
\begin{tabular}{|c|c|c|c|c|}
\hline
   & $n\sqrt{2(1-y_{n,2})}$ & $n\sqrt{2(1-y_{n,3})}$  & $n\sqrt{2(1-y_{n,4})}$ \\ \hline
 $n=150$  & $1.77464$ & $6.0132$ & $9.53661$            \\ \hline
 $n=250$  & $1.10344$ & $5.71202$ & $9.35539$              \\ \hline
 $n=500$  & $1.00403$ & $5.58651$ & $9.27349$                            \\ \hline
\bf{Limit} &$\mathbf{t_{0,1}}=0.903528$ &$\mathbf{t_{0,2}}=5.34057 $ & $\mathbf{t_{0,3}} =9.07889 $    \\ \hline
\end{tabular}}
\end{center}
\begin{center}
\caption[Tabla 8]{\textbf{Case $\gamma=61/5=2(\alpha+2j+1),$} $M=10^6>V$\\ }\label{tabla8}
 $j=3,\ \alpha=-9/10,\ \beta=-9/10,\ \gamma=61/5,\ M_n=\frac{10^6n^2(n-1/2)(n+2)}{n^{\gamma+4}}.$
\end{center}
\end{table}

Finally, we illustrate Theorem \ref{amh} plotting the curves corresponding to the limit functions and to the scaled polynomials $Q_n\left(1-\frac{x^2}{2n^2}\right)$ with  $n=150$ and $n=500.$ In all the figures we have used the same values for the parameters as those ones taken previously in the numerical experiments about the zeros.

\begin{figure}[h]
\centering
\includegraphics[width=3in,height=1.2in]{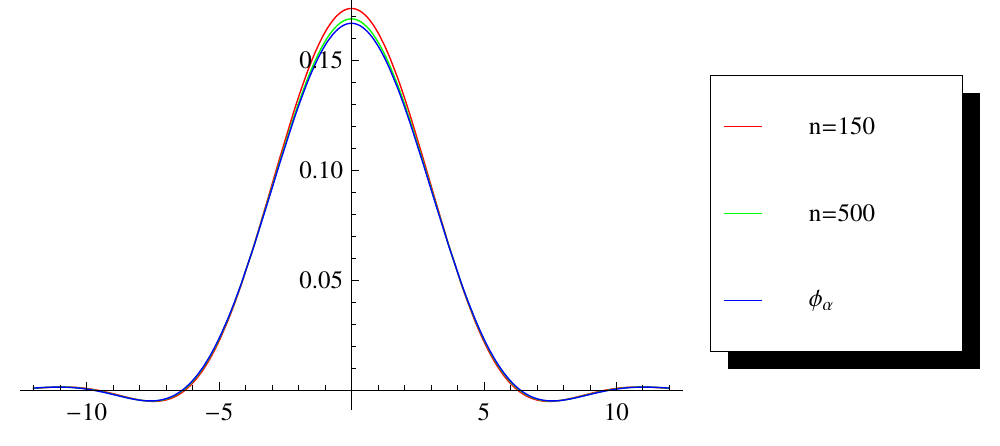}
\caption[Figure 8]{Case $\gamma>2(\alpha+2j+1).$ Limit function and scaled polynomials $Q_n(1-x^2/(2n^2)).$ }\label{greater}
\end{figure}
\begin{figure}[h]
\centering
\includegraphics[width=3in,height=1.2in]{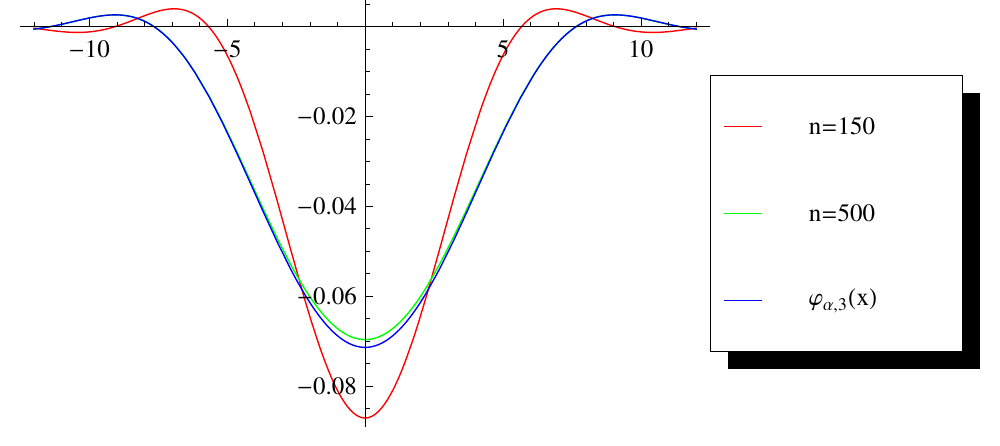}
\caption[Figure 8]{Case $\gamma<2(\alpha+2j+1).$ Limit function and scaled polynomials $Q_n(1-x^2/(2n^2)).$}\label{lesser}
\end{figure}
\begin{figure}[h]
\centering
\includegraphics[width=3in,height=1.2in]{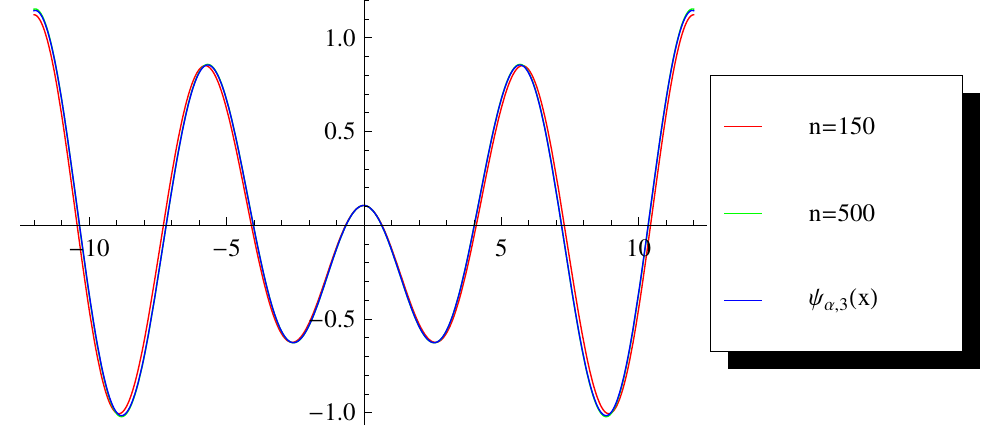}
\caption[Figure 8]{Case $\gamma=2(\alpha+2j+1). $ Limit function and scaled polynomials $Q_n(1-x^2/(2n^2))$ with $M<V.$}\label{equalsmall}
\end{figure}
\begin{figure}[h]
\centering
\includegraphics[width=3in,height=1.2in]{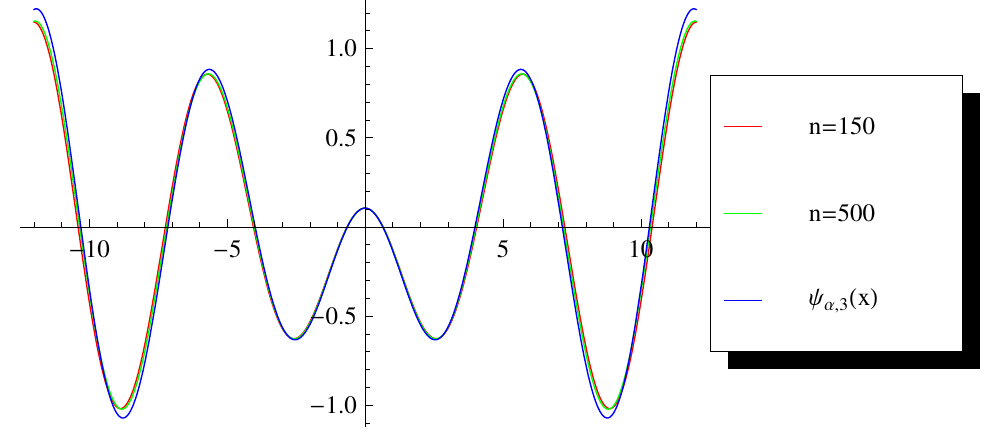}
\caption[Figure 8]{Case $\gamma=2(\alpha+2j+1)$ Limit function and scaled polynomials $Q_n(1-x^2/(2n^2)) $ with $M>V.$}\label{equalbig}
\end{figure}

\noindent  \textbf{Acknowledgements:} We thank the two anonymous referees for their useful suggestions to improve the paper. In special, one  of the referees provided us with some references such as \cite{dls2015} which have been relevant for our research.

\end{document}